\documentclass[fleqn]{mat01}
\usepackage{times,mathtimy,amssymb,latexsym,amsbsy}
\begin{document}

\setcounter{page}{411} \firstpage{411}

\newtheorem{theore}{Theorem}
\renewcommand\thetheore{\arabic{section}.\arabic{theore}}
\newtheorem{theor}[theore]{\bf Theorem}
\newtheorem{rem}[theore]{Remark}
\newtheorem{propo}[theore]{\rm PROPOSITION}
\newtheorem{lem}[theore]{Lemma}
\newtheorem{definit}[theore]{\rm DEFINITION}
\newtheorem{coro}[theore]{\rm COROLLARY}
\newtheorem{exampl}[theore]{Example}
\newtheorem{case}{Case}

\def\corol{\trivlist \item[\hskip \labelsep{COROLLARY.}]}
\def\noteproof{\trivlist \item[\hskip \labelsep{\it Note added in Proof.}]}
\def\proof{\trivlist \item[\hskip \labelsep{\it Proof.}]}
\def\rema{\trivlist \item[\hskip \labelsep{\it Remark.}]}

\def\d{\mbox{\rm d}}
\def\e{\mbox{\rm e}}

\renewcommand{\theequation}{\thesection\arabic{equation}}

\def\A{{\mathcal A}}
\def\C{{\mathbb C}}
\def\F{{\mathbb F}}
\def\G{{\mathcal G}}
\def\H{{\mathbb H}}
\def\Hq{{H(\F_q)}}
\def\N{{\mathbb N}}
\def\P{{\Phi}}
\def\R{{\mathbb R}}
\def\Z{{\mathbb Z}}
\def\a{{\alpha}}
\def\b{{\beta}}
\def\e{{\epsilon}}
\def\p{{\phi}}

\def\PSL{{\operatorname{PSL}}}
\def\PSO{{\operatorname{PSO}}}
\def\PSp{{\operatorname{PSp}}}

\title{On the orders of finite semisimple groups}

\markboth{Shripad M~Garge}{On the orders of finite semisimple groups}

\author{SHRIPAD M~GARGE}

\volume{115}

\mon{November}

\parts{4}

\pubyear{2005}

\address{School of Mathematics, Tata Institute of Fundamental Research, Colaba, Mumbai~400~005, India\\
\noindent E-mail: shripad@math.tifr.res.in}

\Date{MS received 2 April 2005; revised 20 July 2005}

\begin{abstract}
The aim of this paper is to investigate the order coincidences among the
finite semisimple groups and to give a reasoning of such order
coincidences through the transitive actions of compact Lie groups.

It is a theorem of Artin and Tits that a finite simple group is
determined by its order, with the exception of the groups
$(A_3(2), A_2(4))$ and $(B_n(q), C_n(q))$ for $n \geq 3, q$ odd.
We investigate the situation for finite semisimple groups of Lie
type. It turns out that the order of the finite group $H(\F_q)$
for a split semisimple algebraic group $H$ defined over $\F_q$,
does not determine the group $H$ up to isomorphism, but it
determines the field $\F_q$ under some mild conditions. We then
put a group structure on the pairs $(H_1, H_2)$ of split
semisimple groups defined over a fixed field $\F_q$ such that the
orders of the finite groups $H_1(\F_q)$ and $H_2(\F_q)$ are the
same and the groups $H_i$ have no common simple direct factors. We
obtain an explicit set of generators for this abelian,
torsion-free group. We finally show that the order coincidences
for some of these generators can be understood by the inclusions
of transitive actions of compact Lie groups.
\end{abstract}


\maketitle

\section{Introduction}

It is a theorem of Artin and Tits that two finite simple groups of
the same order are isomorphic except for the pairs
\begin{align*}
(\PSL_4(\F_2), \PSL_3(\F_4)) \ \ \, {\rm and} \ \ \, (\PSO_{2n +
1}(\F_q), \PSp_{2n}(\F_q)) \ \ \, {\rm for} \ n \geq 3, q \ {\rm
odd}.
\end{align*}
This theorem was first proved by Emil Artin in 1955 for the finite
simple groups that were known then \cite{A1,A2}. As new finite
simple groups were discovered, Tits \cite{T1,T2,T3,T4,T5,T6}
verified that the above pairs are the only pairs of non-isomorphic
finite simple groups of the same order. One may also look in
\cite{KLST} for an exposition of these proofs.

We investigate in this paper the situation for the groups of
$\F_q$-rational points of a split semisimple algebraic group $H$ defined
over a finite field $\F_q$. Since the orders of the groups $\Hq$ and
$H'(\F_q)$ are the same if $H$ and $H'$ are isogenous and since the
simply connected group is unique in an isogeny class, we concentrate
only on simply connected groups. Since the groups $B_n(\F_q)$ and
$C_n(\F_q)$ have the same order for $n \geq 3$ and for all $q$, we do
not distinguish between them.

This paper is arranged as follows. We state some preliminary
lemmas in \S2 which are used in the proofs of the main theorems.
Section~3 is devoted to determining the field $\F_q$ from the
order of the group $\Hq$. The first natural step in that direction
is to obtain its characteristic. We prove that under certain mild
conditions, the order of the group $\Hq$ determines the
characteristic of $\F_q$ (Theorem~\ref{thm:characteristic}). We
further prove that it eventually determines the base field
(Theorem~\ref{thm:field}).

\setcounter{section}{3} \setcounter{theore}{1}
\begin{theor}[\!] Let $H_1$ and $H_2$
be two split semisimple simply connected algebraic groups defined
over finite fields $\F_{q_1}$ and $\F_{q_2}$ respectively. Let $X$
denote the set $\{8, 9, 2^r, p\}${\rm ,} where $2^r + 1$ is a
prime and $p$ is a prime of the form $2^s \pm 1$. Suppose that{\rm
,} for $i = 1, 2, A_1$ is not a direct factor of $H_i$ whenever
$q_i \in X$ and $B_2$ is not a direct factor of $H_i$ whenever
$q_i = 3$. Then{\rm ,} if $|H_1(\F_{q_1})| = |H_2(\F_{q_2})|${\rm
,} the characteristics of $\F_{q_1}$ and $\F_{q_2}$ are equal.
\end{theor}

Since $|A_1(\F_9)| = |B_2(\F_2)|$, the above theorem is not true
in general. We feel that this is the only counterexample, i.e.,
the conclusion of Theorem~\ref{thm:characteristic} is true without
the hypothesis imposed there except that we must exclude the case
of $H_1 = A_1$ over $\F_9$ and $H_2 = B_2$ over $\F_2$, but we
have not been able to prove it.

\setcounter{section}{3}
\begin{theor}[\!]
Let $H_1$ and $H_2$ be two split semisimple simply connected
algebraic groups defined over finite fields $\F_{q_1}$ and
$\F_{q_2}$ of the same characteristic. Suppose that the order of
the finite groups $H_1(\F_{q_1})$ and $H_2(\F_{q_2})$ are the
same{\rm ,} then $q_1 = q_2$. Moreover the fundamental degrees
(and their multiplicities) of the Weyl groups $W(H_1)$ and
$W(H_2)$ are the same.
\end{theor}

\setcounter{section}{3}
\begin{theor}[\!] Let $H_1$ and $H_2$ be two
split semisimple simply connected algebraic groups defined over a
finite field $\F_q$. If the finite groups $H_1(\F_q)$ and
$H_2(\F_q)$ have the same order then the orders of the groups
$H_1(\F_{q'})$ and $H_2(\F_{q'})$ are the same for any finite
extension $\F_{q'}$ of $\F_q$.
\end{theor}
\setcounter{section}{1}

Thus, the question now boils down to classifying the split
semisimple simply connected groups $H_1, H_2$ defined over a fixed
field $\F_q$ such that the orders of the finite groups $H_1(\F_q)$
and $H_2(\F_q)$ are the same. We first characterise such pairs
where each of the groups $H_1$ and $H_2$ can be written as a
direct product of exactly two simple groups. We find that all such
pairs can be generated by a `nice' set of pairs, which admit a
geometric reason for the order coincidence. We make further
observations regarding the pairs of order coincidence at the end
of \S4.

These observations lead us to a natural question which we answer
in the affirmative in \S5. This question is about describing all
the pairs $(H_1, H_2)$ of groups defined over a fixed finite field
$\F_q$, where the groups $H_i$ have no common simple factor and
the orders of the groups $H_1(\F_q)$ and $H_2(\F_q)$ are the same.
The set of such pairs admits a structure of an abelian,
torsion-free group. We determine an explicit set of generators for
this group (Theorem~\ref{thm:gen}).

Section~6 deals with a geometric reasoning for the `nice' set of pairs
of order coincidence given by transitive action of compact Lie groups.
If $H$ is a compact Lie group and $H_1, H_2$ are connected subgroups of
$H$ such that the natural action of $H_2$ on $H/H_1$ is transitive, then
it can be seen that the split forms of $H \times (H_1 \cap H_2)$ and
$H_1 \times H_2$ have the same number of $\F_q$-rational points for any
finite field $\F_q$. We give such geometric reasoning for the first
three pairs described in Theorem~\ref{thm:gen}. It would be interesting
to know if other pairs also admit such geometric reasoning.

\section{Preliminary lemmas}

\setcounter{equation}{0}
\setcounter{theore}{0}

We state some preliminary lemmas in this section. The first three lemmas
are proved by Artin in his papers \cite{A1,A2}.

Let $\P_n(x)$ be the $n$th cyclotomic polynomial and
\begin{equation*}
\P_n(x, y) = y^{\varphi(n)} \P_n(x/y)
\end{equation*}
be the corresponding homogeneous form. Let $a, b$ be integers which are
relatively prime and which satisfy the inequalities
\begin{equation*}
|a| \geq |b| + 1 \geq 2.
\end{equation*}
Fix a prime $p$ which divides $a^n - b^n$ for some $n$.
Then it is clear that $p$ does not divide any of $a$ and $b$.
Let $f$ be the order of $ab^{-1}$ modulo $p$.
For a natural number $m$, we put ${\rm ord}_p m = \alpha$, where $p^{\alpha}$ is the largest power of $p$ dividing $m$.
We call $p^{{\rm ord}_p m}$ as the $p$-contribution to $m$.

\begin{lem}\hskip -.3pc { $($Lemma $1$ of $[1]$$)$.} \ \ \label{lem:order}
With the above notations{\rm ,} we have the following rules\hbox{{\rm :}}
\leftskip .1pc
(1) If $p$ is odd{\rm ,}
\begin{equation*}
\hskip -1.25pc {\rm ord}_p \P_f(a, b) > 0; \quad {\rm ord}_p
\P_{fp^i}(a, b) = 1 \ {\rm for} \ i \geq 1
\end{equation*}
and in all other cases ${\rm ord}_p \P_n(a, b) = 0$.

Therefore{\rm ,} we have
\begin{align*}
\hskip -1.25pc {\rm ord}_p (a^n - b^n) &= 0,\quad \hbox{if}\quad f \nmid n\\
\hskip -1.25pc {\rm ord}_p (a^n - b^n) &= {\rm ord}_p (a^f - b^f)
+ {\rm ord}_p n,\quad \hbox{otherwise.}
\end{align*}

(2) If $p = 2${\rm ,} then $f = 1$.

\leftskip .1pc
\bigskip
(a) If $\P_1(a, b) = a - b \equiv 0 \pmod
4${\rm ,} then ${\rm ord}_2 \P_{2^i}(a, b) = 1 $ for $i \geq 1$.

(b) If $\P_2(a, b) = a + b \equiv 0 \pmod 4${\rm ,} then ${\rm
ord}_2 \P_{2^i}(a, b) = 1 $ for $i = 0, 2, 3, \dots$.

In all other cases ${\rm ord}_2 \P_n(a, b) = 0$.\vspace{-.9pc}
\end{lem}
\bigskip
\begin{lem}\hskip -.3pc {\rm (\S 4, {\it eqs} (1)--(3) {\it of} \cite{A2}).}\ \ \label{lem:estimate}
Let $\a = (a - 1)(a^2 - 1) \cdots (a^l - 1)$ for some integer $a
\ne 0$ and let $p_1$ be a prime dividing $\a$. Let $P_1$ be the
$p_1$-contribution to $\a${\rm ,} i.e.{\rm ,} $P_1$ be the highest
power of $p_1$ dividing $\a$ and let $q$ be a prime power. We
have\hbox{{\rm :}}

\begin{enumerate}
\renewcommand\labelenumi{{\rm (\arabic{enumi})}}
\leftskip .1pc
\item If $a = \pm q${\rm ,} then $P_1 \leq 2^l (q + 1)^l$.

\item If $a = q^2${\rm ,} then $P_1 \leq 4^l (q + 1)^l$.
\end{enumerate}
\end{lem}

\begin{lem}\hskip -.3pc {\rm ({\it Corollary to Lemma} 2 {\it of} \cite{A1}).}\ \ \label{cor:artin}
If $a > 1$ is an integer and $n > 2$ then there is a prime $p$ which
divides $\P_n(a)$ but no $\P_i(a)$ with $i < n$ unless $n = 6$ and $a =
2$.\vspace{.5pc}
\end{lem}

\begin{lem}\label{lem:induction}
If the inequality $q^n \geq \alpha (q + 1)${\rm ,} where $\alpha$
is a fixed positive real number{\rm ,} holds for a pair of
positive integers $(q_1, n_1)${\rm ,} then it holds for all $(q_2,
n_2)$ satisfying $q_2 \geq q_1$ and $n_2 \geq n_1$.\vspace{.5pc}
\end{lem}

\begin{proof}
This is clear. \hfill $\Box$
\end{proof}

\begin{lem}\label{lem:herbrand}
Let $H$ be a semisimple algebraic group defined over a finite field
$\F_q$. If $\widetilde{H}$ denotes a connected cover of $H${\rm ,} then
$|\widetilde{H}(\F_q)| = |\Hq|$.
\end{lem}

\begin{proof}
Let $\widetilde{H}$ be a connected cover of $H$.
We have an exact sequence
$0 \rightarrow A \rightarrow \widetilde{H} \rightarrow H \rightarrow 1$
where $A$ is a finite abelian group.
From this sequence, we get the following exact sequence of Galois cohomology sets
\begin{equation*}
0 \rightarrow H^0(\F_q, A) \rightarrow \widetilde{H}(\F_q)
\rightarrow H(\F_q) \rightarrow H^1(\F_q, A) \rightarrow H^1(\F_q,
\widetilde{H}).
\end{equation*}
By Lang's Theorem (Corollary to theorem~1 of \cite{L}), $H^1(\F_q,
\widetilde{H}) = 0$. Since all the sets in the above sequence are
finite, we have
\begin{equation*}
|H^0(\F_q, A)| \cdot |H(\F_q)| = |\widetilde{H}(\F_q)| \cdot
|H^1(\F_q, A)|.
\end{equation*}
Since the Galois group, ${\rm Gal}(\overline{\F}_q/\F_q)$, is
procyclic and $A$ is a finite Galois-module, its Herbrand quotient
is $1$, i.e., $|H^0(\F_q, A)| = |H^1(\F_q, A)|$. It follows that
$|\widetilde{H}(\F_q)| = |H(\F_q)|$.

\hfill $\Box$\vspace{-1.2pc}
\end{proof}

\section{Determining the finite field}

\setcounter{equation}{0}
\setcounter{theore}{0}

The first natural step in determining the field $\F_q$ is to
determine its characteristic. Observe that if we have two
semisimple groups $H_1$ and $H_2$ defined over finite fields
$\F\!_{p_1^{r_1}}$ and $\F\!_{p_2^{r_2}}$ respectively, such that
$|H_1(\F\!_{p_1^{r_1}})| = |H_2(\F\!_{p_2^{r_2}})|$ and $p_1 \ne
p_2$, then either $p_1$ fails to give the largest contribution to
the order of $H_1(\F\!_{p_1^{r_1}})$ or $p_2$ fails to give the
largest contribution to the order of $H_2(\F\!_{p_2^{r_2}})$.
Therefore we would like to obtain a description of the split
semisimple algebraic groups $H$ defined over $\F\!_{p^r}$ such
that the $p$-contribution to the order of the group
$H(\F\!_{p^r})$ is not the largest. These groups are the only
possible obstructions for determining the characteristic of the
base field. Since we limit ourselves to the case of simply
connected groups only, every semisimple group considered in this
paper is a direct product of (simply connected) simple algebraic
groups. Hence we need to describe simple algebraic groups $H$
defined over $\F\!_{p^r}$ with the property that $p$ does not
contribute the largest to the order of $H(\F\!_{p^r})$.

We remark that the main tool in the proof of the following
proposition is Lemma~\ref{lem:order} which is proved by Artin in
\cite{A1}. Our proof of the following proposition is very much on
the lines of Artin's proof of Theorem $1$ in \cite{A2}. However,
our result is for $\Hq$, the groups of $\F_q$-rational points of a
simple algebraic group $H$ defined over $\F_q$ whereas Artin
proved the result for finite groups that are simple. The groups
$\Hq$ that we consider here, are not always simple, because of the
presence of (finite) center. Moreover, we get the counterexamples
$A_1(\F_9)$, $A_1(\F_p)$ for a Fermat prime $p$, and $B_2(\F_3)$
which do not figure in Artin's list of counterexamples described
in Theorem~1 of \cite{A2}.

\begin{propo}\label{propn:counter}$\left.\right.$\vspace{.5pc}

\noindent Let $H$ be a split simple algebraic group defined over a
finite field $\F_q$ of characteristic $p$. If the $p$-contribution
to the order of the finite group $\Hq$ is not the largest prime
power dividing the order{\rm ,} then the group $\Hq$ is\hbox{{\rm
:}}\vspace{-.3pc}
\begin{enumerate}
\renewcommand\labelenumi{{\rm (\arabic{enumi})}}
\leftskip .1pc
\item $A_1(\F_q)$ for $q \in \{8, 9, 2^r, p\}$ where $2^r + 1$ is a
Fermat prime and $p$ is a prime of the type $2^s \pm 1$ or
\item $B_2(\F_3)$.\vspace{-.3pc}
\end{enumerate}
Moreover in all these cases{\rm ,} the $p$-contribution is the
second largest prime power dividing the order of the group $\Hq$.
\end{propo}

We call the groups, $A_1(\F_q)$ and $B_2(\F_3)$, described above,
as {\it counterexamples} in the remaining part of this paper.

\begin{proof}
We first recall the orders of the finite groups $\Hq$ where $H$ is
a split simple algebraic group defined over a finite field $\F_q$
(see \S2.9 of \cite{C}).
\begin{align*}
|A_n(\F_q)| &= q^{n(n+1)/2}(q^2 - 1)(q^3 - 1) \cdots (q^{n + 1} - 1), \quad n \geq 1, \\[.15pc]
|B_n(\F_q)| &= q^{n^2}(q^2 - 1)(q^4 - 1) \cdots (q^{2n} - 1), \quad n \geq 2, \\[.15pc]
|D_n(\F_q)| &= q^{n(n-1)}(q^2 - 1)(q^4 - 1) \cdots (q^{2n - 2} - 1)(q^n - 1), \quad n \geq 4,\\[.15pc]
|G_2(\F_q)| &= q^6 (q^2 - 1)(q^6 - 1), \\[.15pc]
|F_4(\F_q)| &= q^{24} (q^2 - 1)(q^6 - 1)(q^8 - 1)(q^{12} - 1), \\[.15pc]
|E_6(\F_q)| &= q^{36} (q^2 - 1)(q^5 - 1)(q^6 - 1)(q^8 - 1)(q^9 - 1)(q^{12} - 1), \\[.15pc]
|E_7(\F_q)| &= q^{63}(q^2 - 1)(q^6 - 1)(q^8 - 1)(q^{10} - 1)(q^{12} - 1)(q^{14} - 1)\\
&\quad \times (q^{18} - 1), \\[.15pc]
|E_8(\F_q)| &= q^{120} (q^2 - 1)(q^8 - 1)(q^{12} - 1)(q^{14} - 1)(q^{18} - 1)(q^{20} - 1)\\
 &\quad \times (q^{24} - 1)(q^{30} - 1).
\end{align*}
Now, let $H$ be one of the finite simple groups listed above and let $p_1$ be a prime dividing the order of the finite group $\Hq$ such that $p_1 \nmid q$.
We use Lemma~\ref{lem:estimate} to estimate $P_1$, the $p_1$-contribution to the order of $\Hq$.
Depending on the type of $\Hq$, we put the following values of $a$ and $l$ in Lemma~\ref{lem:estimate}.
\begin{center}
\begin{tabular}{||c||c|c|c|c|c|c|c||}
\hline
$H$ & $A_n$ & $B_n, D_n$ & $G_2$ & $F_4$ & $E_6$ & $E_7$ & $E_8$ \\\hline
$a$ & $q$ & $q^2$ & $q^2$ & $q^2$ & $q$ & $q^2$ & $q^2$ \\\hline
$l$ & $n + 1$ & $n$ & $3$ & $6$ & $12$ & $9$ & $15$ \\\hline
\end{tabular}
\end{center}

Now, suppose that the $p$-contribution to the order of the group $\Hq$ is not the largest, i.e., the power of $q$ that appears in the formula for $|\Hq|$ is smaller than $P_1$ for some prime $p_1 \nmid q$.
Then, depending on the type of the group, we get following inequalities from Lemma~\ref{lem:estimate}:
\begin{alignat*}{3}
A_n &: P_1 \leq 2^{n + 1}(q + 1)^{n + 1} & &\implies q^{n/2} < 2(q + 1), \\[.15pc]
B_n &: P_1 \leq 4^n (q + 1)^n & &\implies q^n < 4(q + 1),\\[.15pc]
D_n &: P_1 \leq 4^n (q + 1)^n & &\implies q^{n-1} < 4(q + 1),\\[.15pc]
G_2 &: P_1 \leq 4^3 (q + 1)^3 & &\implies q^2 < 4(q + 1),\\[.15pc]
F_4 &: P_1 \leq 4^6 (q + 1)^6 & &\implies q^4 < 4(q + 1),\\[.15pc]
E_6 &: P_1 \leq 2^{12} (q + 1)^{12} & &\implies q^3 < 2(q + 1),\\[.15pc]
E_7 &: P_1 \leq 4^9 (q + 1)^9 & &\implies q^7 < 4(q + 1),\\[.15pc]
E_8 &: P_1 \leq 4^{15} (q + 1)^{15} & &\implies q^8 < 4(q + 1).
\end{alignat*}
In all the cases where the above inequalities of the type $q^m < \alpha (q + 1)$ do not hold, we get that $p$ contributes the largest to the order of $\Hq$.
Observe that the last four inequalities, i.e., the inequalities corresponding to the groups $F_4$, $E_6$, $E_7$ and $E_8$ do not hold for $q = 2$ and hence by
Lemma~\ref{lem:induction} they do not hold for any $q \geq 2$.
Thus, for $H = F_4, E_6, E_7$ and $E_8$, the $p$-contribution to $|\Hq|$ is always the largest prime power dividing $|\Hq|$.

Similarly we obtain the following table of the pairs of positive
integers $(q, n)$ where the remaining inequalities fail. Then using
Lemma~\ref{lem:induction}, we know that for all $(q', n')$ with $q' \geq
q$ and $n' \geq n$, the contribution of the characteristic to the order
of the finite group $\Hq$ is the largest. Therefore, we are left with
the cases for $(q', n')$ such that $q' < q$ or $n' < n$, which are to be
checked. The adjoining table shows the groups $\Hq$ which are to be
checked.
\begin{center}
\begin{tabular}{|r|ll|}\hline
$A_n$ & $q = 2$, & $n \geq 6$\\
 & $q = 3, 4, 5$, &$n \geq 4$\\
 & $q \geq 7$, & $n \geq 3$\\
$B_n$ & $q = 2$, & $n \geq 4$\\
 & $q = 3, 4$, & $n \geq 3$\\
 & $q \geq 5$, & $n \geq 2$\\
$D_n$ & $q = 2$, & $n \geq 5$\\
 & $q \geq 3$, & $n \geq 4$\\
$G_2$ & $q \geq 5$ &\\\hline
\end{tabular}
\hskip7mm
\begin{tabular}{|l|}
\hline
$A_3(\F_2), A_3(\F_3), A_3(\F_4), A_3(\F_5),$ \\
$A_4(\F_2), A_5(\F_2),$ \\
$A_1(\F_q), A_2(\F_q) \hskip2mm \forall q,$ \\
$B_2(\F_2), B_2(\F_3), B_2(\F_4),$ \\
$B_3(\F_2),$ \\
 \\
$D_4(\F_2),$ \\
 \\
$G_2(\F_2), G_2(\F_3), G_2(\F_4)$. \\\hline
\end{tabular}
\end{center}
In all the cases other than $A_1(\F_q)$ and $A_2(\F_q)$, we can do
straightforward calculations and check that $p$ contributes the largest
to the order of every group except for $B_2(\F_3)$. In the case of
$B_2(\F_3)$, the prime $3$ indeed fails to give the largest
contribution, however it gives the second largest contribution to the
order of the group. The cases of $A_1$ and $A_2$ over a general finite
field $\F_q$ are done in a different way.

We first deal with the case of the group $A_2(\F_q)$.
Recall that
\begin{equation*}
|A_2(\F_q)| = q^3 (q^2 - 1)(q^3 - 1) = q^3 (q^2 + q + 1)(q + 1)(q - 1)^2.
\end{equation*}
Let $p_1 \nmid q$ be a prime dividing the order of $A_2(\F_q)$ and
let $P_1$ be the contribution of $p_1$ to $|A_2(\F_q)|$. Let $f$
denote the order of $q$ modulo $p_1$. If $f \neq 1$ then it is
clear that the $p_1$-contribution to the order of $A_2(\F_q)$ is
not more than $q^3$. If $f = 1$ and $p_1 \ne 2$ or $3$, then by
Lemma~\ref{lem:order}, $P_1$ divides $(q - 1)^2$ which is less
than $q^3$. If $f = 1$ and $p_1 = 2$ or $3$, then $P_1$ divides
either $3(q - 1)^2$, $2(q - 1)^2$ or $4(q + 1)$. Thus, if $q^3$ is
not the largest prime power dividing $|A_2(\F_q)|$, then $q^3 <
P_1$ for some prime $p_1 \ne p$ and hence we have
\begin{equation*}
q^3 < 3(q - 1)^2, \quad 2(q - 1)^2 \quad {\rm or}\quad 4(q + 1).
\end{equation*}
Again as above, we observe that none of the above inequalities are
satisfied by $q \geq 3$, and then we check that the $2$-contribution to
the order of $A_2(\F_2)$ is the largest one.

Now, for the group $A_1(\F_q)$, we observe that for any prime $p_1 \nmid
q$, the $p_1$-contribution to the order of $A_1(\F_q)$ divides $q^2 - 1
= (q + 1)(q - 1)$. We make two cases here depending on $q$ being odd or
even.

If $q$ is odd, both $q + 1$ and $q - 1$ are even. The $2$-contribution
to one of the numbers $q + 1$ and $q - 1$ is $2$, and the other number
then must be a power of $2$ if $q$ is not the largest prime power
dividing $|A_1(\F_q)|$. If $q + 1$ is a power of $2$, then $q$ is
necessarily a prime of the form $q = p = 2^s - 1$, a Mersenne prime.
However, if $q - 1$ is a power of $2$ then the only possibilities for
$q$ are that $q$ is a Fermat prime, $q = p = 2^s + 1$, or $q = 9$.

If $q$ is even, $q - 1$ and $q + 1$ are both odd and hence they do
not have any common prime factor. If $q$ is not the largest prime
power dividing $|A_1(\F_q)|$, then the largest prime power
dividing $|A_1(\F_q)|$ must be $q + 1$. Let $q = 2^r$ and $P_1 =
p_1^s = 2^r + 1$. Here $p_1$ is odd, and hence by
Lemma~\ref{cor:artin}, if $s >2$, there is a prime divisor of
$p_1^s - 1$ which does not divide $p_1 - 1$, a contradiction. If
$s = 2$, $2^r = p_1^2 - 1$. Then both $p_1 \pm 1$ are powers of
two and hence we obtain that $p_1 = 3$ and $q = 2^3 = 8$. If $s =
1$, $p_1 = 2^r + 1$, a Fermat\break prime. \hfill $\Box$
\end{proof}

Thus, if $\Hq$ is not one of the counterexamples described in the
above proposition, then the characteristic of $\F_q$ contributes
the largest to the order of $\Hq$. Since every (simply connected)
semisimple algebraic group is a direct product of (simply
connected) simple algebraic groups, we get that whenever a finite
semisimple group $\Hq$ does not have any of the above
counterexamples as direct factors, then the characteristic of
$\F_q$ contributes the largest to the order of $\Hq$.

\begin{theor}[\!]\label{thm:characteristic}
Let $H_1$ and $H_2$ be two split semisimple simply connected algebraic
groups defined over finite fields $\F_{q_1}$ and $\F_{q_2}$
respectively. Let $X$ denote the set $\{8, 9, 2^r, p\}$ where $2^r + 1$
is a Fermat prime and $p$ is a prime of the type $2^s \pm 1$. Suppose
that for $i = 1, 2, A_1$ is not one of the direct factors of $H_i$
whenever $q_i \in X$ and $B_2$ is not a direct factor of $H_i$ whenever
$q_i = 3$. Then{\rm ,} if $|H_1(\F_{q_1})| = |H_2(\F_{q_2})|${\rm ,} the
characteristics of $\F_{q_1}$ and $\F_{q_2}$ are equal.
\end{theor}

\begin{proof}
This is clear.\hfill $\Box$
\end{proof}

Now, we come to the main theorem of this section.
Recall that if $H$ is a split semisimple algebraic group of rank $n$ defined over a finite field $\F_q$, then the order of $\Hq$ is given by the formula,
\begin{equation*}
|\Hq| = q^N (q^{d_1} - 1) (q^{d_2} - 1) \cdots (q^{d_n} - 1),
\end{equation*}
where $d_1, d_2, \dots, d_n$ are the degrees of the basic
invariants of $W(H)$, the Weyl group of $H$ and $N = \sum_i (d_i -
1)$ (\S2.9 of \cite{C}). Now onwards, we call $d_i$ as the {\it
degrees} of $W(H)$, the Weyl group of $H$. Observe that for every
split simple algebraic group $H$, the integer $2$ always occurs as
a degree of $W(H)$ with multiplicity one (\S3.7 of \cite{H}).
Therefore the multiplicity of the integer $2$ among the degrees of
$W(H)$ determines the number of simple direct factors of the group
$H$. We remark here that the degrees $d_i$ of $W(H)$ may appear
with multiplicities.

\begin{theor}[\!]\label{thm:field}
Let $H_1$ and $H_2$ be two split semisimple simply connected
algebraic groups defined over finite fields $\F_{q_1}$ and
$\F_{q_2}$ of the same characteristic. Suppose that the order of
the finite groups $H_1(\F_{q_1})$ and $H_2(\F_{q_2})$ are the
same{\rm ,} then $q_1 = q_2$. Moreover the fundamental degrees
{\rm (}and the multiplicities{\rm )} of the Weyl groups $W(H_1)$
and $W(H_2)$ are the same.
\end{theor}

\begin{proof}
Let $p$ be the characteristic of the fields $\F_{q_1}$ and $\F_{q_2}$, and let $q_1 = p^{t_1}, q_2 = p^{t_2}$.
Let the orders of the finite groups $H_1(\F_{q_1})$ and $H_2(\F_{q_2})$ be given by
\begin{align*}
|H_1(\F_{q_1})| &= (q_1)^r~(q_1^{r_1} - 1) (q_1^{r_2} - 1) \cdots (q_1^{r_n} - 1)
\\[.3pc]
&= (p^{t_1})^r ((p^{t_1})^{r_1} - 1) ((p^{t_1})^{r_2} - 1) \cdots ((p^{t_1})^{r_n} - 1)\\[.3pc]
|H_2(\F_{q_2})| &= (q_2)^s (q_2^{s_1} - 1)~(q_2^{s_2} - 1) \cdots (q_2^{s_m} - 1) \\[.3pc]
&= (p^{t_2})^s ((p^{t_2})^{s_1} - 1)~ ((p^{t_2})^{s_2} - 1) \cdots
((p^{t_2})^{s_m} - 1)
\end{align*}
As remarked above, the integers $r_i$ and $s_j$ are the respective
degrees of the Weyl groups $W(H_1)$ and $W(H_2)$. Moreover the
rank of the group $H_1$ is $n$ and that of $H_2$ is $m$. Further,
we have
\begin{equation*}
r = \sum_{i = 1}^n (r_i - 1) \quad {\rm and} \quad s = \sum_{j = 1}^m (s_j - 1).
\end{equation*}
Since $|H_1(\F_{q_1})| = |H_2(\F_{q_2})|$, we have that
\begin{equation*}
t_1r = t_2s
\end{equation*}
and
\begin{equation}\label{eqn:1}
\prod_{i = 1}^n ((p^{t_1})^{r_i} - 1) = \prod_{j = 1}^m
((p^{t_2})^{s_j} - 1).
\end{equation}
Assume that $r_1 \leq r_2 \leq \cdots \leq r_n$ and $s_1 \leq s_2
\leq \cdots \leq s_m$. We treat both the products in
eq.~\eqref{eqn:1} as polynomials in $p$ and factor them into the
cyclotomic polynomials in $p$.

Let us assume for the time being that $p \ne 2$, so that we can
apply Lemma \ref{cor:artin} to conclude that the cyclotomic
polynomials appearing on both sides of eq.~\eqref{eqn:1} are the
same with the same multiplicities. Observe that on the left-hand
side (LHS) the highest order cyclotomic polynomial is
$\P_{t_1r_n}(p)$ whereas such a polynomial on the right-hand side
(RHS) is $\P_{t_2s_m}(p)$. Since the cyclotomic polynomials
appearing on both sides are the same, we have that $t_1r_n =
t_2s_m$. Thus, the polynomial $p^{t_1r_n} - 1$, which is same as
the polynomial $p^{t_2s_m} - 1$, can be cancelled from both sides
of eq.~\eqref{eqn:1}. Continuing in this way we get that
$t_1r_{n-k} = t_2s_{m-k}$ for all $k$. This implies in particular
that $m = n$. Further,
\begin{equation*}
t_1r = t_2s \implies \sum t_1r_i - t_1n = \sum t_2s_j - t_2n.
\end{equation*}
But, by the above observation, this gives us that $t_1n = t_2n$
and hence $t_1 = t_2$, i.e., $q_1 = q_2$. Thus, the fields
$\F_{q_1}$ and $\F_{q_2}$ are isomorphic.

Now, it also follows that $r_i = s_i$ for all $i$, i.e., the degrees of
the corresponding Weyl groups are the same.

Now, let $p = 2$. So, we have the equation
\begin{equation}\label{eqn:2}
(2^{t_1})^r \prod_{i = 1}^n ((2^{t_1})^{r_i} - 1) = (2^{t_2})^s
\prod_{j = 1}^m ((2^{t_2})^{s_j} - 1).
\end{equation}
The only possible obstruction to the desired result in this case
comes from $\P_6(2) = 3$ and $\P_2(2) =3$. Moreover, if $\P_6(2)$
divides the LHS of the equation but not the RHS, then it is clear
that $(2^3 - 1)(2^2 - 1)^2$ divides the RHS with the same power as
that of $2^6 - 1$ in the LHS. Other than these polynomials, all
the factors of type $2^l - 1$ occur on both sides with the same
multiplicities.

Since $s_j > 1$ for all $j$, $t_2 = 1$, i.e., $q_2 = 2$ and the possible values for $q_1$ are $2$, $2^2$ and $2^3$, since $t_1$ divides $6$.
We prove the result in only one case, $q_1 = 2$, as other cases can be handled by a similar reasoning.
If $q_1 = 2$, eq.~\eqref{eqn:2} becomes
\begin{equation*}
2^r \prod_{i = 1}^n (2^{r_i} - 1) = 2^s \prod_{j = 1}^m (2^{s_j} -
1).
\end{equation*}
But then $r = s$ and hence $\sum_i (r_i - 1) = \sum_j (s_j - 1)$.
Now, the term $(2^6 - 1)$ contributes $5$ to $r$ whereas the term
$(2^3 - 1)(2^2 - 1)^2$ contributes only $4$ to $s$. As other
factors are same on both the sides, this is a contradiction. Hence
the factor $(2^6 - 1)$ in LHS of eq.~\eqref{eqn:2} must be
adjusted by the same factor in the RHS. Then we get that $m = n$
and $r_i = s_i$ for all $i$. \hfill $\Box$
\end{proof}

\begin{theor}[\!]\label{thm:weyl}
Let $H_1$ and $H_2$ be two split semisimple simply connected algebraic
groups defined over a finite field $\F_q$. If the orders of the finite
groups $H_1(\F_q)$ and $H_2(\F_q)$ are same then the orders of
$H_1(\F_{q'})$ and $H_2(\F_{q'})$ are the same for any finite extension
$\F_{q'}$ of $\F_q$.
\end{theor}

\begin{proof}
Let $H_1$, $H_2$ be split semisimple algebraic groups defined over
$\F_q$. By Theorem~\ref{thm:field}, we have that if $|H_1(\F_q)| =
|H_2(\F_q)|$ then the degrees of the Weyl groups $W(H_1)$ and
$W(H_2)$ are the same with the same multiplicities. Then the
formulae for the orders of the groups $H_1(\F_q)$ and $H_2(\F_q)$
are the same as polynomials in $q$. Hence the orders of the groups
$H_1(\F_{q'})$ and $H_2(\F_{q'})$ are the same for any finite
extension $\F_{q'}$ of $\F_q$. \hfill $\Box$\vspace{-.3pc}
\end{proof}

\section{Order coincidences}

\setcounter{equation}{0}
\setcounter{theore}{0}

We fix a finite field $\F_q$ and all algebraic groups considered in this
section are assumed to be defined over $\F_q$.

In this section, we concentrate on the pairs of split semisimple
groups $(H_1, H_2)$, such that the orders of the groups
$H_1(\F_q)$ and $H_2(\F_q)$ are the same. We want to characterise
all possible pairs of order coincidence $(H_1, H_2)$ and to
understand the reason behind the coincidence of these orders. This
section and the next one are devoted towards characterising all
pairs of order coincidence and we discuss a geometric reasoning of
these order coincidences in the last section.

We know by Theorem~\ref{thm:field}, that for such a pair $(H_1,
H_2)$, the degrees of the corresponding Weyl groups, $W(H_1)$ and
$W(H_2)$, must be the same with the same multiplicities. For a
Weyl group $W$, we denote the collection of degrees of $W$ by
$d(W)$. We now make the following easy observations which follow
from the basic theory of the Weyl groups \cite{H}.

\begin{rem}\label{rem:degrees}
{\rm Let $H_1$ and $H_2$ be two split semisimple algebraic groups over a finite field $\F_q$ such that the groups $H_1(\F_q)$ and $H_2(\F_q)$ have the same order.
Then we have:\vspace{-.3pc}
\begin{enumerate}
\renewcommand\labelenumi{{\rm (\arabic{enumi})}}
\leftskip .1pc
\item The rank of the group $H_1$ is the same as the rank of $H_2$.
\item The number of direct simple factors of the groups $H_1$ and $H_2$ is the same.
\item If one of the groups, say $H_1$, is simple, then so is $H_2$ and in that case $H_1$ is isomorphic to
$H_2$.\vspace{-.3pc}
\end{enumerate}
(We remind the reader once again that we do not distinguish between the groups of type $B_n$ and $C_n$.)}
\end{rem}

The next natural step would be to look at the order coincidences
in the case of groups each having two simple factors. We
characterise such pairs in the following theorem.

\begin{theor}[\!]\label{thm:pairs}
Let $H_1$ and $H_2$ be split semisimple simply connected algebraic
groups each being a direct product of exactly two simple algebraic
groups. Assume that $H_1$ and $H_2$ do not have any common simple direct
factor. Then the pairs $(H_1, H_2)$ such that $|H_1(\F_q)| =
|H_2(\F_q)|$ are exhausted by the following list\hbox{{\rm :}} \vspace{-.3pc}
\begin{enumerate}
\renewcommand\labelenumi{{\rm (\arabic{enumi})}}
\leftskip .1pc
\item $(A_{2n - 2}B_n, A_{2n - 1}B_{n - 1})$ for $n \geq 2${\rm ,} with the convention that $B_1 = A_1${\rm ,}
\item $(A_{n - 2}D_n, A_{n - 1}B_{n - 1})$ for $n \geq 4${\rm ,}
\item $(B_{n - 1}D_{2n}, B_{2n - 1}B_n)$ for $n \geq 2${\rm ,} with the convention that $B_1 = A_1${\rm ,}
\item $(A_1A_5, A_4G_2)${\rm ,}
\item $(A_1B_3, B_2G_2)${\rm ,}
\item $(A_1D_6, B_5G_2)${\rm ,}
\item $(A_2B_3, A_3G_2)$ and
\item $(B_3^2, D_4G_2)$.\vspace{-.5pc}
\end{enumerate}
\end{theor}

\begin{proof}
Let $H_1 = H_{1, 1} \times H_{1, 2}$ and $H_2 = H_{2, 1} \times H_{2,
2}$ where $H_{i, j}$ are split simple algebraic groups. We denote
$W(H_i)$ by $W_i$ and $W(H_{i, j})$ by $W_{i, j}$.

Since the orders of the groups $|H_1(\F_q)|$ and $|H_2(\F_q)|$ are the
same, by Theorem~\ref{thm:field}, the degrees of the Weyl groups $W_1$
and $W_2$ are the same with the same multiplicities. Moreover for $i =
1, 2$, we have $W_i = W_{i, 1} \times W_{i, 2}$.

Let $n$ be the maximum of the degrees of $W_1$. Then it is the
largest of the maxima of the degrees of $W_{1, j}$ for $j = 1, 2$.
Suppose that $n$ is the maximum degree of $W_{1, 1}$. Then
depending on $n$, we have the following choices for the group
$H_{1, 1}$:
\begin{center}
\begin{tabular}{|c|l||c|l|}
\hline
$n = 2$ &$A_1,$ & $n = 4$ & $A_3, B_2,$ \\
$n = 6$ &$A_5, B_3, D_4, G_2,$ & $n = 12$ & $A_{11}, B_6, D_7, F_4, E_6,$ \\
$n = 18$ &$A_{17}, B_9, D_{10}, E_7,$ & $n = 30$ & $A_{29}, B_{15}, D_{16}, E_8,$ \\
$n$ odd &$A_{n - 1},$ & $\begin{matrix} n = 2m  \\
n \not\in \{2, 4, 6, 12, 18, 30\} \end{matrix}$ &$A_{2m - 1}, B_m, D_{m + 1},$ \\ \hline
\end{tabular}
\end{center}

The general philosophy of the proof is as follows:

Once we fix $n$, we have a finite set of choices for $H_{1, 1}$
and $H_{2, 1}$. Then we fix one choice each for $H_{1, 1}$ and
$H_{2, 1}$, and compare the degrees of $W_{1, 1}$ and $W_{2, 1}$.
Since $H_{1, 1} \neq H_{2, 1}$ the collections $d(W_{1, 1})$ and
$d(W_{2, 1})$ are different. The degrees of $W_{1, 1}$ that do not
occur in $d(W_{2, 1})$ must occur in the collection $d(W_{2, 2})$
and similarly the degrees of $W_{2, 1}$ that do not occur in
$d(W_{1, 1})$ must occur in the collection $d(W_{1, 2})$. Moreover
the degrees of $W_{1, 2}$ and $W_{2, 2}$ are bounded above by $n$.
This gives us further finitely many choices for the groups $H_{1,
2}$ and $H_{2, 2}$. For these choices, we simply verify the
equality of the collections $d(W_1)$ and $d(W_2)$. If the
collections are equal, we get a coincidence of orders $(H_1,
H_2)$.

As a sample, we do the case of $n = 4$ to illustrate the above philosophy.

Let us assume that $H_{1, 1} = A_3$ and $H_{2, 1} = B_2$. Then we have
\begin{equation*}
d(W_{1, 1}) = \{2, 3, 4\} \subseteq d(W_1) \quad {\rm and} \quad
d(W_{2, 1}) = \{2, 4\} \subseteq d(W_2).
\end{equation*}
Thus, $3$ is a degree of $W_{2, 2}$ and the maximum of the degrees
of $W_{2, 2}$ is less than or equal to $4$, hence $H_{2, 2} =
A_2$. Then, since $d(W_1) = d(W_2)$, the only possibility for the
collection $d(W_{1, 2})$ is $\{2\}$ and we get $H_{1, 2} = A_2$.
This gives us the order coincidence $(A_1A_3, A_2B_2)$.

\hfill $\Box$\vspace{-.5pc}
\end{proof}

Observe that in the above theorem, we have three infinite families
of pairs,\vspace{-.5pc}
\begin{align*}
&(A_{2n - 2}B_n, A_{2n - 1}B_{n - 1}),  (A_{n - 2}D_n, A_{n -
1}B_{n - 1})
\end{align*}
and\vspace{-.5pc}
\begin{align*}
&(B_{n - 1}D_{2n}, B_{2n - 1}B_n).
\end{align*}
If we consider the following pairs given by the first two infinite families:
\begin{align*}
(H_1, H_2) &= (A_{2n - 2}B_n, A_{2n - 1}B_{n - 1})
\end{align*}
and\vspace{-.5pc}
\begin{align*}
(H_3, H_4) &= (A_{2n - 2}D_{2n}, A_{2n - 1}B_{2n - 1}),
\end{align*}
then\vspace{-.5pc}
\begin{equation*}
(H_1H_4, H_2H_3) = (A_{2n - 2}A_{2n - 1}B_nB_{2n - 1}, A_{2n - 1}A_{2n - 2}B_{n - 1}D_{2n}).
\end{equation*}
This implies that $(B_{2n - 1}B_n, B_{n - 1}D_{2n})$ is also a
pair of order coincidence and this is precisely our third infinite
family! Thus, the third infinite family of order coincidences can
be obtained from the first two infinite families.

Similarly if we consider
\begin{equation*}
(H_1, H_2) = (A_2D_4, A_3B_3) \quad {\rm and} \quad (H_3, H_4) =
(A_2B_3, A_3G_2),
\end{equation*}
then we get the pair $(B_3^2, D_4G_2)$ from the pair $(H_2H_3, H_1H_4)$.

Similarly we observe that
\begin{align*}
(A_1B_3, B_2G_2) \ \hbox{can be obtained from} \ (A_1A_3, A_2B_2) \ \hbox{and} \ (A_2B_3, A_3G_2),\\[.2pc]
(A_1A_5, A_4G_2) \ \hbox{can be obtained from} \ (A_4B_3, A_5B_2) \ \hbox{and} \ (A_1B_3, B_2G_2),\\[.2pc]
(A_1D_6, B_5G_2) \ \hbox{can be obtained from} \ (A_4D_6, A_5B_5)
\ \hbox{and} \ (A_1A_5, A_4G_2).
\end{align*}

We record our observation as a remark below.

\begin{rem}\label{rem:pairs}
{\rm All the pairs of order coincidence described in Theorem~\ref{thm:pairs} can be obtained from the following pairs:
\begin{enumerate}
\renewcommand\labelenumi{{\rm (\arabic{enumi})}}
\leftskip .1pc
\item $(A_{2n - 2}B_n, A_{2n - 1}B_{n - 1})$ for $n \geq 2$, with the convention that $B_1 = A_1$,\vspace{.1pc}
\item $(A_{n - 2}D_n, A_{n - 1}B_{n - 1})$ for $n \geq 4$, and\vspace{.1pc}
\item $(A_2B_3, A_3G_2)$.\vspace{-.4pc}
\end{enumerate}}
\end{rem}

These pairs are quite special, in the sense that they admit a
geometric reasoning for the coincidence of orders. We describe it
in the last section.

If we do not restrict ourselves to the groups having exactly two
simple factors, then we also find the following pairs $(H_1, H_2)$
involving other exceptional groups:
\begin{align*}
(A_1B_4B_6, B_2B_5F_4), (A_4G_2A_8B_6, A_3A_6B_5E_6), (A_1B_7B_9,
B_2B_8E_7),
\end{align*}
and\vspace{-.5pc}
\begin{align*}
(A_1B_4B_7B_{10}B_{12}B_{15}, ~B_3B_5B_8B_{11}B_{14}E_8).
\end{align*}
One now asks a natural question whether these four pairs, together with
the pairs described in Remark~\ref{rem:pairs}, generate all possible
pairs of order coincidence. We make this question more precise in the
next section and answer it in the affirmative.

\section{On a group structure on pairs of groups of equal order}

\setcounter{equation}{0}
\setcounter{theore}{0}

Fix a finite field $\F_q$. Let $\A$ be the set of ordered pairs $(H_1, H_2)$
where $H_1$ and $H_2$ are split semisimple algebraic groups defined over the
field $\F_q$ such that the orders of the finite groups $H_1(\F_q)$ and
$H_2(\F_q)$ are the same. We define an equivalence relation on $\A$ by saying
that an element $(H_1, H_2) \in \A$ is related to $(H_1', H_2') \in \A$, denoted
by $(H_1, ~H_2) \sim (H_1', ~H_2')$, if and only if there exist two split semisimple
algebraic groups $H$ and $K$ defined over $\F_q$ such that
\begin{equation*}
H_1' \times K = H_1 \times H \quad {\rm and} \quad H_2' \times K = H_2 \times H.
\end{equation*}
It can be checked that $\sim$ is an equivalence relation. We
denote the set of equivalence classes in $\A$ given by $\sim,
\A/{\sim}$, by $\G$ and the equivalence class of an element $(H_1,
H_2) \in \A$ is denoted by $[(H_1, H_2)]$. This set $\G$ describes
all pairs of order coincidence $(H_1, H_2)$ where the split
semisimple (simply connected) groups $H_i$ do not have any common
direct simple factor.

We put a binary operation on $\G$ given by
\begin{equation*}
[(H_1, H_2)] \circ [(H_1', H_2')] = [(H_1 \times H_1', H_2 \times
H_2')].
\end{equation*}
It is easy to see that the above operation is a well-defined
modulo, the equivalence that we have introduced. The set $\G$ is
obviously closed under $\circ$ which is an associative operation.
The equivalence class $[(H, H)]$ acts as the identity and $[(H_1,
H_2)]^{-1} = [(H_2, H_1)]$. Thus $\G$ is an abelian, torsion-free
group. Since the first two infinite families described in
Remark~\ref{rem:pairs} are independent, the group $\G$ is not
finitely generated.

Let $\G'$ be the subgroup of $\G$ generated by following elements.
\begin{enumerate}
\renewcommand\labelenumi{{\rm (\arabic{enumi})}}
\leftskip .1pc
\item ${\sf B}_n = [(A_{2n - 2}B_n, A_{2n - 1}B_{n - 1})]$, for $n \geq 2$, with the convention that $B_1 = A_1$, \vspace{.1pc}
\item ${\sf D}_n = [(A_{n - 2}D_n, A_{n - 1}B_{n - 1})]$, for $n \geq 4$, \vspace{.1pc}
\item ${\sf G}_2 = [(A_2B_3, ~A_3G_2)]$, \vspace{.1pc}
\item ${\sf F}_4 = [(A_1B_4B_6, ~B_2B_5F_4)]$, \vspace{.1pc}
\item ${\sf E}_6 = [(A_4G_2A_8B_6, ~A_3A_6B_5E_6)]$, \vspace{.1pc}
\item ${\sf E}_7 = [(A_1B_7B_9, ~B_2B_8E_7)]$, \vspace{.1pc}
\item ${\sf E}_8 = [(A_1B_4B_7B_{10}B_{12}B_{15}, ~B_3B_5B_8B_{11}B_{14}E_8)]$.
\end{enumerate}
(For a group such as $B_n$, we use ${\sf B}_n$ to denote a pair
$(H_1, ~H_2)$ in which $B_n$ appears as a group of the largest
degree.)

\begin{lem}\label{lem:gen}
Let $n$ be a positive integer. Let $H_1$ and $H_2$ be split{\rm ,}
simply connected{\rm ,} simple algebraic groups such that the Weyl
groups $W(H_1)$ and $W(H_2)$ have the same highest degree and it
is equal to $n$. Then there is an element in $\G'$ which can be
represented as the equivalence class of a pair $(K_1, K_2)$ such
that for $i = 1, 2${\rm ,} $H_i$ is one of the simple factors of
$K_i$ and for any other simple factor $H_i'$ of $K_i$ the highest
degree of $W(H_i')$ is less than $n$.
\end{lem}

\begin{proof}
We prove this lemma by explicit calculations. If $n$ is odd or $n
= 2$, there is nothing to prove as there is only one group, $A_{n
-1}$, with $n$ as the highest degree.

For $n = 4$, the groups $A_3$ and $B_2$ are the only groups with
$4$ as the highest degree and ${\sf B}_2 = [(A_2B_2, ~A_3B_1)]$ is
an element of the group $\G'$ where $A_3$ and $B_2$ appear as
factors on either sides and all other simple groups that appear
have highest degree less than $4$.

If $n = 2m$ for $m > 2$ and  $m \not \in \{3, 6, 9, 15\}$, then $A_{2m - 1}$, $B_{m}$ and $D_{m + 1}$ are the only groups with $n$ as the highest degree.
Consider following elements of $\G'$:
\begin{align*}
&{\sf B}_m = [(A_{2m - 2}B_m, A_{2m - 1}B_{m - 1})],\\
&{\sf D}_{m + 1} = [(A_{m - 1}D_{m + 1}, A_mB_m)]
\end{align*}
and\vspace{-.5pc}
\begin{align*}
&{\sf D}_{m + 1} \circ {\sf B}_m = [(A_{m - 1}A_{2m - 2}D_{m + 1},
~A_mA_{2m - 1}B_{m - 1})].
\end{align*}
The element ${\sf B}_m$ contains the simple groups $A_{2m - 1}$
and $B_m$ on its either sides and other simple groups appearing in
${\sf B}_m$ have highest degree less than $2m$. Similarly the
elements ${\sf D}_{m + 1}$ and ${\sf D}_{m + 1} \circ {\sf B}_m$
are the required elements of $\G'$ for the pairs $\{D_{m + 1},
B_m\}$ and $\{A_{2m - 1}, D_{m + 1}\}$.

Now, we consider the cases when $n \in \{6, 12, 18, 30 \}$. These
cases involve exceptional groups.

For $n = 6$, the groups $A_5, B_3, D_4$ and $G_2$ are the only groups
with $6$ as the highest degree. We have following elements of $\G'$ for
the corresponding pairs.
\begin{align*}
&{\sf B}_3 = [(A_4B_3, A_5B_2)] \ \ \hbox{{\rm for the pair}} \ \ \{B_3, A_5\},\\
&{\sf D}_4 = [(A_2D_4, A_3B_3)] \ \ \hbox{{\rm for the pair}} \ \ \{D_4, B_3\},\\
&{\sf G}_2 = [(A_2B_3, A_3G_2)] \ \ \hbox{{\rm for the pair}} \ \ \{B_3, G_2\},\\
&{\sf D}_4 \circ {\sf G}_2 = [(A_2^2D_4, A_3^2G_2)] \ \ \hbox{{\rm for the pair}} \ \ \{D_4, G_2\},\\
&{\sf B}_3 \circ {\sf D}_4 = [(A_2A_4D_4, A_3A_5B_2)] \ \ \hbox{{\rm for the pair}} \ \ \{D_4, A_5\},\\
&{\sf G}_2 \circ {\sf B}_3^{-1} = [(A_2A_5B_2, A_3A_4G_2)] \ \
\hbox{{\rm for the pair}} \ \ \{A_5, G_2\}.
\end{align*}

In the same way, we give the following elements of the group $\G'$
for all possible simple groups having highest degree $12, 18$ and
$30$.\vspace{.5pc}
\begin{center}
\begin{tabular}{llll@{}}\hline
 & & &\\[-.8pc]
Element of $\G'$ &Pair &Element of $\G'$ &Pair\\\hline
 & & &\\[-.5pc]
\multicolumn{4}{c}{$n = 12$}\\
 & & &\\[-.5pc]
${\sf B}_6 = [(A_{10}B_6, A_{11}B_5)]$ & $\{B_6, A_{11}\}$, & ${\sf D}_7 = [(A_5D_7, A_6B_6)]$ & $\{D_7, B_6\}$, \\[.25pc]
${\sf F}_4 = [(A_1B_4B_6, B_2B_5F_4)]$ & $\{B_6, F_4\}$, & ${\sf E}_6 = [(A_4G_2A_8B_6,$ &$\{B_6, E_6\}$\\
 & &$\qquad \ A_3A_6B_5E_6)]$ &\\[.25pc]
${\sf B}_6 \circ {\sf D}_7$ & $\{D_7, A_{11}\}$, & ${\sf D}_7 \circ {\sf F}_4$ & $\{D_7, F_4\}$, \\[.25pc]
${\sf B}_6^{-1} \circ {\sf F}_4$ & $\{A_{11}, F_4\}$, & ${\sf D}_7 \circ {\sf E}_6$ & $\{D_7, E_6\}$, \\[.25pc]
${\sf B}_6^{-1} \circ {\sf E}_6$ & $\{A_{11}, E_6\}$, & ${\sf F}_4^{-1} \circ {\sf E}_6$ & $\{F_4, E_6\}$.\\[.5pc]\hline
 & & &\\[-.5pc]
\multicolumn{4}{c}{$n = 18$}\\
 & & &\\[-.5pc]
${\sf B}_9 = [(A_{16}B_9, A_{17}B_8)]$ & $\{B_9, A_{17}\}$, & ${\sf D}_{10} = [(A_8D_{10}, A_9B_9)]$ & $\{D_{10}, B_9\}$,\\[.25pc]
${\sf E}_7 = [(A_1B_7B_9, B_2B_8E_7)]$ & $\{B_9, E_7\}$, & ${\sf B}_9^{-1} \circ {\sf E}_7$ & $\{A_{17}, E_7\}$, \\[.25pc]
${\sf B}_9 \circ {\sf D}_{10}$ & $\{D_{10}, A_{17}\}$, & ${\sf D}_{10} \circ {\sf E}_7$ & $\{D_{10}, E_7\}$. \\[.5pc]\hline
 & & &\\[-.5pc]
\multicolumn{4}{c}{$n = 30$}\\
 & & &\\[-.5pc]
${\sf B}_{15} = [(A_{28}B_{15}, A_{29}B_{14})]$ & $\{B_{15}, A_{29}\}$, & ${\sf B}_{15}^{-1} \circ {\sf E}_8$ & $\{A_{29}, E_8\}$, \\[.25pc]
${\sf D}_{16} = [(A_{14}D_{16}, A_{15}B_{15})]$ & $\{D_{16}, B_{15}\},$ & ${\sf B}_{15} \circ {\sf D}_{16}$ & $\{D_{16}, A_{29}\},$ \\[.25pc]
${\sf E}_8 = [(A_1B_4B_7B_{10}B_{12}B_{15},$ & $\{B_{15}, E_8\},$ & ${\sf D}_{16} \circ {\sf E}_8$ & $\{D_{16}, E_8\}.$\\
$\qquad \ B_3B_5B_8B_{11}B_{14}E_8)]$ & & &\\[.3pc]\hline
\end{tabular}
\end{center}
This completes the proof of the lemma. \hfill $\Box$
\end{proof}

\begin{theor}[\!]\label{thm:gen}
The groups $\G$ and $\G'$ are the same. In other words{\rm ,} the group $\G$ is
generated by the following elements\hbox{{\rm :}}
\begin{enumerate}
\renewcommand\labelenumi{{\rm (\arabic{enumi})}}
\leftskip .1pc
\item $[(A_{2n - 2}B_n, A_{2n - 1}B_{n - 1})]$ for $n \geq 2${\rm ,}
with the convention that $B_1 = A_1${\rm ,}\vspace{.1pc}
\item $[(A_{n - 2}D_n, A_{n - 1}B_{n - 1})]$ for $n \geq 4${\rm ,}\vspace{.1pc}
\item $[(A_2B_3, A_3G_2)]${\rm ,} \vspace{.1pc}
\item $[(A_1B_4B_6, B_2B_5F_4)]${\rm ,} \vspace{.1pc}
\item $[(A_4G_2A_8B_6, A_3A_6B_5E_6)]${\rm ,} \vspace{.1pc}
\item $[(A_1B_7B_9, B_2B_8E_7)]${\rm ,} \vspace{.1pc}
\item $[(A_1B_4B_7B_{10}B_{12}B_{15}, B_3B_5B_8B_{11}B_{14}E_8)]$.
\end{enumerate}
\end{theor}

\begin{proof}
Let $[(H_1, H_2)] \in \G$. By Theorem~\ref{thm:field}, the
fundamental degrees of the Weyl groups $W(H_1)$ and $W(H_2)$ are
the same with the same multiplicities. Let $n$ be the highest
degree of $W(H_1)$ which is the same as the highest degree of
$W(H_2)$. For $i = 1, 2$, let $K_i$ be one of the simple factors
of $H_i$ such that $n$ is the highest degree of $W(K_i)$. Then by
previous lemma, there exists an element $[(H_1', H_2')] \in \G'$
such that $K_i$ are the simple factors of $H_i'$ and the other
simple factors of $H_i'$ have highest degree less than $n$. Thus,
the element $[(H_1H_2', H_2H_1')]$ is an element of the group $\G$
and the multiplicity of $K_i$ on either side of this element is
now reduced by $1$. This way, we cancel all the simple factors
having $n$ as the highest degree and then the result is obtained
by\break induction.\hfill $\Box$
\end{proof}

\section{Transitive actions of compact Lie groups}

\setcounter{equation}{0}
\setcounter{theore}{0}

Here we explain how a transitive action of compact Lie groups is
related to the coincidence of orders. The exposition is based on
Chapter~$2$, page 121 of \cite{GOV}.

Suppose $H$ is a compact simply connected Lie group acting
transitively on a compact manifold $X = H/H_1$ with $H_1$
connected. Suppose that $H_2$ is a closed connected Lie subgroup
of $H$ and that the action of $H$ on $X$ when restricted to $H_2$
remains transitive. Then $X = H/H_1 = H_2/(H_1 \cap H_2)$. By
looking at the homotopy exact sequence for the fibration $1
\rightarrow H' \rightarrow H \rightarrow H/H' \rightarrow 1$ for
any closed subgroup $H'$ of $H$,
\begin{equation*}
\pi_1(H') \rightarrow \pi_1(H) \rightarrow \pi_1(H/H') \rightarrow
\pi_0(H'),
\end{equation*}
we find that $H/H'$ is simply connected if and only if $H'$ is
connected. Therefore $X = H/H_1$ is simply connected and hence if
$X = H_2/(H_1 \cap H_2)$ with $H_2$ simply connected, $H_1 \cap
H_2$ is connected.

We now assume that there is an analogue of the action of $H$ on $X
= H/H_1$ over finite fields, which we now take to be all defined
over $\F_q$. By Lang's theorem (Corollary to Theorem~1 of
\cite{L}) if $H_1$ is connected then
\begin{equation*}
|(H/H_1)(\F_q)| = \frac{|H(\F_q)|}{|H_1(\F_q)|}.
\end{equation*}
Therefore for the equality of spaces $H/H_1$ and $H_2/(H_1 \cap
H_2)$, with $H_1, H_2, H_1 \cap H_2$ connected, we find that
\begin{equation*}
|H(\F_q)| \cdot |(H_1 \cap H_2)(\F_q)| = |H_1(\F_q)| \cdot |H_2(\F_q)|.
\end{equation*}
Thus transitive action of compact Lie groups gives rise to
coincidence of orders of finite semisimple groups.

We call an ordered $3$-tuple $(H, H_1, H_2)$, as discussed above,
a triple of {\em inclusion of transitive actions}. We first
classify all such triples of inclusion of transitive actions and
explain the geometric reasoning behind the order coincidence for
the first three pairs described in Theorem~\ref{thm:gen}. We note
some observations.

\begin{rem}\label{rem:triple}
{\rm Let ($H, H_1, H_2$) be a triple of inclusion of transitive actions,
where $H$, $H_1$ and $H_2$ are compact Lie groups such that $H_1$ is a
subgroup of $H$ and the natural action of $H_1$ on $H/H_2$ is
transitive. Then
\begin{enumerate}
\renewcommand\labelenumi{{\rm (\arabic{enumi})}}
\leftskip .1pc
\item $H = H_1H_2$ (Lemma~4.1, page~138 of \cite{GOV}) and
\item either $H_1$ or $H_2$ has the same maximal exponent as the maximal
exponent of the group $H$ (Corollary~2, page~143 of
\cite{GOV}).\vspace{-.4pc}
\end{enumerate}
(We recall that a natural number $a$ is an exponent of a compact Lie group
$H$ if and only if $a + 1$ is a degree of the Weyl group of the split form
of $H$.)}
\end{rem}

Therefore to classify the inclusions among the transitive actions,
equivalently to determine the triples $(H, H_1, H_2)$ of inclusion of
transitive actions, it would be desirable to classify the subgroups of a
given Lie group of the maximal exponent. We restrict ourselves to the
case when $H$ is a simple Lie group.

\begin{theor}[\cite{O}] \label{thm:exponent}
Let $H$ be a connected simple compact Lie group and $H_1$ be a
compact Lie subgroup of $H$ of maximal exponent. Then{\rm ,} the
pairs $H_1 \subseteq H$ are exhausted by the following
list\hbox{{\rm :}}
\begin{align*}
&Sp_n \subset SU_{2n} \quad (n > 1), \quad G_2 \subset SO_7, \quad SO_{2n - 1} \subset SO_{2n} \quad (n > 3),\\
&\quad\ {\rm Spin}_7 \subset SO_8, \quad G_2 \subset SO_8, \quad
F_4 \subset E_6.
\end{align*}
\end{theor}

Observe that the subgroup $H_1 \subset H$ is automatically a
simple group. Now, we classify the triples $(H, H_1, H_2)$, of
inclusion of transitive actions, where $H$ is a simple Lie group.

\begin{theor}[\cite{O}] \label{thm:triple}
The triples $(H, H_1, H_2)$ of inclusion of transitive actions where the
group $H$ is simple are the following ones\hbox{{\rm :}}
\begin{center}
\begin{tabular}{|c|c|c|c||c|c|c|c|}\hline
$H$ & $H_1$ & $H_2$ & $H_1 \cap H_2$ & $H$ & $H_1$ & $H_2$ & $H_1 \cap H_2$ \\\hline
$\begin{matrix} SU_{2n} \\ (n \geq 2) \end{matrix}$ & $Sp_n$ & $SU_{2n - 1}$ & $Sp_{n - 1}$ & $\begin{matrix} SO_{4n} \\ (n \geq 2) \end{matrix}$ & $SO_{4n - 1}$ & $Sp_n$ & $Sp_{n - 1}$ \\ \hline
$SO_7$ & $G_2$ & $\begin{matrix} SO_6 \\ SO_5 \end{matrix}$ & $\begin{matrix} SU_3 \\ SU_2 \end{matrix}$ & $SO_{16}$ & $SO_{15}$ & ${\rm Spin}_9$ & ${\rm Spin}_7$ \\\hline
$\begin{matrix} SO_{2n} \\ (n \geq 4) \end{matrix}$ & $SO_{2n - 1}$ & $SU_n$ & $SU_{n - 1}$ & $SO_8$ & ${\rm Spin}_7$ & $\begin{matrix} SO_7 \\ SO_6 \\ SO_5 \end{matrix}$ & $\begin{matrix} G_2 \\ SU_3 \\ SU_2 \end{matrix}$ \\\hline
\end{tabular}
\end{center}
\end{theor}

Observe that the exponents of the groups $H \times (H_1 \cap H_2)$
and $H_1 \times H_2$ are the same in all the above cases. Hence
$(H \times (H_1 \cap H_2), ~H_1 \times H_2)$ is a pair of order
coincidence for us. The pairs described in Remark~\ref{rem:pairs}
occur in the above descriptions. We can, in fact, give an explicit
description of the inclusion among the transitive actions
corresponding to the pairs given in Remark~\ref{rem:pairs}.

The groups $O_n$, $U_n$ and $Sp_n$ act on the spaces $\R^n$, $\C^n$ and
$\H^n$, respectively, in a natural way. By restricting this action to
the corresponding spheres, we get that the groups $O_n$, $U_n$ and
$Sp_n$ act transitively on the spheres $S^{n - 1}$, $S^{2n - 1}$ and
$S^{4n - 1}$, respectively. By fixing a point in each of the spheres, we
get the corresponding stabilizers as $O_{n - 1} \subset O_n$, $U_{n - 1}
\subset U_n$ and $Sp_{n - 1} \subset Sp_n$.

By treating the space $\C^n = \R^{2n}$, we obtain an inclusion of
transitive actions $U_n \subset O_{2n}$, with both the groups acting
transitively on $S^{2n - 1}$. Since $S^{2n - 1}$ is connected, the
actions of $SU_n \subset U_n$ and $SO_{2n} \subset O_{2n}$ on $S^{2n -
1}$ remain transitive. Thus we get an inclusion of actions $SU_n \subset
SO_{2n}$ and the corresponding stabilisers are $SU_{n - 1}\subset SO_{2n
- 1}$. Thus, we get a triple $(SO_{2n}, SU_n, SO_{2n - 1})$ or
equivalently we get a pair of order coincidence as $(D_nA_{n - 2}, ~A_{n
- 1}B_n)$.

Similarly, by treating $\H^n$ as $\C^{2n}$ and repeating the above
arguments, we get the inclusion of transitive actions $Sp_n \subset
SU_{2n}$, acting on the sphere $S^{4n - 1}$, with $Sp_{n - 1} \subset
SU_{2n - 1}$ as the corresponding stabilisers. This gives us the triple
$(SU_{2n}, Sp_n, SU_{2n - 1})$ and the pair of order coincidence $(A_{2n
- 1}B_{n - 1}, B_nA_{2n - 2})$.

Thus, we get the two infinite families described in
Remark~\ref{rem:pairs}. The remaining pair of order coincidence,
$(A_1B_3, B_2G_2)$, can be obtained in a similar way by
considering the natural inclusion $G_2 \subset SO_7$. These groups
act transitively on the sphere $S^6$ and the corresponding
stabilisers are $SU_3 \subset SO_6$. We observe that the split
form of $SO_6$ is isomorphic to $SL_4$, and therefore the triple
$(SO_7, G_2, SO_6)$ gives us $(A_2B_3, A_3G_2)$ as the
corresponding pair of order coincidence.

\begin{rem}
{\rm It would be interesting to know if the pairs $(4)$ to $(7)$ of
Theorem~\ref{thm:gen} involving exceptional groups are also obtained
in this geometric way.}
\end{rem}

\section*{Acknowledgements}

The author thanks Prof. Dipendra Prasad for suggesting the
question and for many fruitful discussions. The author also thanks
Dr Maneesh Thakur for many stimulating discussions. He also thanks
TIFR, Mumbai and the Bhaskaracharya Pratishthana, Pune, where a
part of this work was done.

\end{document}